\documentclass[11pt]{amsart}
\usepackage{amssymb}
\usepackage{amsmath}
\usepackage{amscd}
\usepackage[all]{xypic}
\usepackage{graphicx}

\allowdisplaybreaks[1]

\newtheorem{thm}{Theorem}[section]
\newtheorem{lem}[thm]{Lemma}
\newtheorem{cor}[thm]{Corollary}
\newtheorem{prop}[thm]{Proposition}

\newtheorem{tadothm}{Theorem}[section]

\theoremstyle{definition}
\newtheorem{rem}[thm]{Remark}
\newtheorem{defn}[thm]{Definition}

\newtheorem{ex}[thm]{Example}

\theoremstyle{remark}

\numberwithin{equation}{section}

\def\R{{\mathbb R}}
\def\Z{{\mathbb Z}}
\def\C{{\mathbb C}}

\def\Hom{\text{\rm Hom}}

\newcommand{\ft}[3]{{f_{#1}\otimes f_{#2}\otimes f_{#3}}}

\begin{document}
\title[Pointed harmonic volumes of hyperelliptic curves]
{The pointed harmonic volumes of hyperelliptic curves
with Weierstrass base points}

\author[Yuuki Tadokoro]{Yuuki Tadokoro}
\address{Department of Mathematical Sciences,
University of Tokyo, 3-8-1 Komaba, Meguro, Tokyo 153-8914, Japan}
\email{tado\char`\@ms.u-tokyo.ac.jp}


\begin{abstract}
We give a explicit computation of
the pointed harmonic volumes of hyperelliptic curves
with Weierstrass base points,
which are paraphrased into a combinatorial formula.
\end{abstract}

\maketitle
\section{Introduction}
Let $X$ be a compact Riemann surface
of genus $g\geq 2$
and $p$ a point on $X$.
By Pulte \cite{P}, the pointed harmonic volume of $(X,p)$ was defined
to be the homomorphism $K\otimes H\to \R/{\Z}$,
using Harris' method for the harmonic volume of $X$ \cite{H-1}.
Here, we denote by $H=H^1(X;\Z)$ the
first cohomology group of $X$ 
and $K$ the kernel of the intersection pairing
$H\otimes H\to \Z$.
In this paper, we compute the pointed harmonic
volume of any hyperelliptic curve $C$ with
any Weierstrass point $p$.
In theorem \ref{a pointed harmonic volume},
we compute that of some special hyperelliptic curve
$C_0$ with Weierstrass points in an analytic way,
by the explicit computation of
Chen's iterated integrals \cite{C}.
Using Proposition \ref{key fact},
we can compute the the pointed
harmonic volumes of all the hyperelliptic curves with
Weierstrass base points from
those of $C_0$.
These results are paraphrased from
a combinatorial viewpoint as follows.
Let $\{P_{j}\}_{j=0,1,\ldots,2g+1}$ denote
the set of Weierstrass points on $C$,
and fix a Weierstrass point $P_{\nu}$,
$0\leq \nu\leq 2g+1$.
A certain homomorphism
$\kappa_{\nu}\colon K\otimes H\to \dfrac{1}{2}\Z/{\Z}=\{0,1/{2}\}$
is defined in \S \ref{formula},
which depends on the choice of $P_\nu$.
\begin{tadothm}
For any hyperelliptic curve $C$ and $A\in K\otimes H$, we have
$$
I_{P_\nu}(A)\equiv\kappa_{\nu}(A)
\ \mathrm{mod}\ \Z.
$$
\end{tadothm}
The author \cite{T} computed the harmonic
volumes of hyperelliptic curves.
But the computation of the pointed ones
of $(X,p)$ is more complicated
than that of $X$.
For any hyperelliptic curve $C$,
it is tedious to compute $I_p$
in the case $p\in C\setminus\{P_j\}_{j=0,1,\ldots,2g+1}$.
But we have $I_p\equiv 0$ or $1/{2}$ $\mathrm{mod}\ \Z$
in the case $p\in \{P_j\}_{j=0,1,\ldots,2g+1}$.
It has been still unknown
which elements of $K\otimes H$ and Weierstrass points $p$
have nontrivial $I_p$ or not.
In this paper, we compute them completely.

As an application of the pointed harmonic volume of $(X,p)$,
Pulte proved the pointed Torelli theorem \cite{P}.
We denote by $\pi_1(X,p)$ the fundamental group of
$X$ at the base point $p\in X$ and
$J_p$ the augmentation ideal of the group ring $\Z\pi_1(X,p)$.
\begin{thm}$($the pointed Torelli theorem \cite{P}$)$\ \\
Suppose that $X$ and $Y$ are compact Riemann surfaces and
that $p\in X$ and $q\in Y$.
With the exception of two points $p$ in $X$,
if there is a ring isomorphism
$$\Z\pi_1(X,p)/{J_p^3}\to \Z\pi_1(Y,q)/{J_q^3}$$
which preserves the mixed Hodge structure,
then there is a biholomorphism
$\varphi\colon X\to Y$ such that $\varphi(p)=q$.
\end{thm}
If $X$ is generic (e.g. $X$ is hyperelliptic),
then there are no exceptional points.
The pointed harmonic volumes determine the choice of
the base points.
In the proof of this theorem,
the classical Torelli theorem follows
from the preservation of the mixed Hodge structure
and we obtain the biholomorphism $X \cong Y$.
When we choose the base points,
the pointed harmonic volume plays
an important role.
Theorem \ref{a combinatorial formula}
also tells the choice of Weierstrass base points
on $C$.

Now we describe the contents of this paper
briefly.
In \S \ref{The pointed harmonic volume},
we define the pointed harmonic volume
of $(X,p)$, using Chen's iterated integrals \cite{C}.
In \S \ref{Hyperelliptic curves},
we give a basis of the first homology
group $H_1(C;\Z/{2\Z})$ of the hyperelliptic curve $C$.
In \S \ref{Pointed harmonic volumes of hyperelliptic curves
and the moduli space of compact Riemann surfaces},
we prove $I_{P_\nu}\in H^0(\Delta^1_g; \Hom(K\otimes H,\Z/{2\Z}))$.
In \S \ref{a hyperelliptic curve},
the pointed harmonic volume of
some special hyperelliptic curve $C_0$
with Weierstrass base points
is computed in an analytic way.
This result can be extended to all the hyperelliptic curves with
Weierstrass base points
and interpreted from a combinatorial viewpoint.
In \S \ref{formula},
we obtain a simple combinatorial formula
of the pointed harmonic volume of $(C,P_{\nu})$.

\noindent
{\bf Acknowledgments.}
The author is grateful to Nariya Kawazumi
for valuable advice
and reading the manuscript.
He also thanks Masahiko Yoshinaga for
reading the manuscript.
This work is partially supported by 21st
Century COE program  (University of Tokyo) by
the Ministry of Education, Culture, Sports, Science and Technology.


\section{The pointed harmonic volume}
\label{The pointed harmonic volume}
We recall the definition of the pointed harmonic volume
of a pointed Riemann surface $(X,p)$.
Here $X$ is a compact Riemann surface
of genus $g\geq 2$
and $p$ a point on $X$.
We identify
the first integral homology group $H_1(X;\mathbb{Z})$
of $X$ with the first integral cohomology group
by Poincar\'e duality,
and denote it by $H$.
For closed $1$-forms $\omega_{1,i}$ and $\omega_{2,i}$,
$i=1,2,\ldots,m$, on $X$ such that
$\displaystyle \int_{X}\sum_{i=1}^{m}
\omega_{1,i}\wedge\omega_{2,i}=0$,
we obtain the $1$-form $\eta$ such that
$d\eta=\sum_{i=1}^{m}\omega_{1,i}\wedge \omega_{2,i}$ and
$\displaystyle \int_{X}\eta\wedge\ast\alpha=0$
for any closed $1$-form $\alpha$ on $X$.
Here, $\ast$ is the Hodge star operator
which depends only on the complex structure
and not the choice of Hermitian metric.
We identify $H$
with the space of all the
real harmonic $1$-forms on $X$
with integral periods by the Hodge theorem.
We denote by $K$ the kernel of the
intersection pairing 
$(\ ,\ )\colon H\otimes H\to \Z$.
\begin{defn}\mbox{(The pointed harmonic volume \cite{P})}\\
For $\sum_{i=1}^{m}a_i\otimes b_i\in K$
and $c\in H$,
the pointed harmonic volume is
defined to be
$$I_{p}{\Biggl(}{\biggl(}\sum_{i=1}^{m}a_{i}\otimes b_{i}{\biggr)}\otimes c{\Biggr)}:=\sum_{i=1}^{m}\int_{\gamma}a_{i}b_{i}-\int_{\gamma}\eta
 \quad \mathrm{mod} \ \mathbb{Z}.$$
Here $\eta$ is the $1$-form on $X$ which
is associated to 
$\sum_{i=1}^{m}a_i\otimes b_i$
in the way stated above and
$\gamma\colon [0,1]\to X$ is a loop in $X$ at the base point $p$
whose homology class is equal to $c$.
The integral $\displaystyle \int_{\gamma}a_ib_i$
is Chen's iterated integral \cite{C}, that is,
$\displaystyle \int_{\gamma}a_ib_i
=\int_{0\leq t_1\leq t_2\leq 1}f_i(t_1)g_i(t_2)dt_1dt_2$
for $\gamma^{\ast}a_i=f_i(t)dt$ and $\gamma^{\ast}b_i=g_i(t)dt$.
Here $t$ is the coordinate in the unit interval $[0,1]$.
See Chen \cite{C} for iterated integrals
and Harris \cite{H-1}, Pulte \cite{P}
for the (pointed) harmonic volume.
\end{defn}
\begin{rem}\label{invarinant}
By the definition of $I_{p}$,
we have
$I_p((\sum_{i=1}^{m}a_{i}\otimes b_{i})\otimes c)
\equiv -I_p((\sum_{i=1}^{m}b_{i}\otimes a_{i})\otimes c)
\ \textnormal{mod}\ \Z$.
\end{rem}

\section{Hyperelliptic curves}
\label{Hyperelliptic curves}
Let $C$ be a hyperelliptic curve
and $\Z_2$ the field $\Z/{2\Z}$.
In this section,
we explain the first homology group of $C$
with $\Z_2$-coefficient.

We define the hyperelliptic curve $C$
as follows.
It is the compactification of the plane curve
in the $(z,w)$ plane $\C^2$
$$w^2=\prod_{i=0}^{2g+1} (z-p_i),$$
where $p_0, p_1,\ldots, p_{2g+1}$ are some
distinct points on $\C$.
It admits the hyperelliptic involution
given by $\iota:(z,w) \mapsto (z,-w)$.
Let $\pi$ be the $2$-sheeted covering 
$C \to \C P^1, (z,w) \mapsto z$, branched over $2g+2$
branch points $\{p_i\}_{i=0,1,\cdots ,2g+1}$
and $P_i \in C$ a ramification point such that $\pi(P_i)=p_i$.
It is known that $\{P_i\}_{i=0,1,\ldots,2g+1}$
is just the set of
all the Weierstrass points on any hyperelliptic curve $C$.

For points $p_i$ and $p_j$,
we denote by $p_ip_j$ a simple path
joining $p_i$ and $p_j$.
We draw simple paths $p_0p_1,p_1p_2,\ldots, p_{2g}p_{2g+1}$
and $p_{2g+1}p_0$ such that
all the $2g+2$ arcs do not intersect
except for endpoints of them.
We take a disk $D\subset \C P^1$
whose boundary is
$\big(\bigcup_{j=0}^{2g}p_jp_{j+1}\big)\bigcup p_{2g+1}p_{0}$
(Figure \ref{disk}, $g=2$).
We picture two copies of $\C P^1$ as above
and call them $\Omega_0$ and $\Omega_1$.
We make crosscuts along
$p_{2k}p_{2k+1}, k=0,1,\ldots,g$
and construct the hyperelliptic curve $C$
by joining 
every $p_{2k}p_{2k+1}$ on $\Omega_0$
to the corresponding one on $\Omega_1$
for $k=0,1,\ldots,g$.
See 102-103 in \cite{F-K} for example.
We may consider $\Omega_i\subset C$ for $i=0,1$.
\begin{figure}[htbp]
\begin{center}
\begin{minipage}{.46\linewidth}
\scalebox{0.65}[0.65]{
\unitlength 0.1in
\begin{picture}( 42.5000, 23.1200)( 13.5000,-26.4200)
%
\special{pn 8}%
\special{sh 0.600}%
\special{ar 3800 600 42 42  0.0000000 6.2831853}%
%
\special{pn 8}%
\special{sh 0.600}%
\special{ar 2600 600 42 42  0.0000000 6.2831853}%
%
\special{pn 8}%
\special{sh 0.600}%
\special{ar 1400 1400 42 42  0.0000000 6.2831853}%
%
\special{pn 8}%
\special{sh 0.600}%
\special{ar 3000 2000 42 42  0.0000000 6.2831853}%
%
\special{pn 8}%
\special{sh 0.600}%
\special{ar 4800 1200 42 42  0.0000000 6.2831853}%
%
\special{pn 8}%
\special{sh 0.600}%
\special{ar 4200 2600 42 42  0.0000000 6.2831853}%
%
\special{pn 13}%
\special{pa 2600 600}%
\special{pa 3800 600}%
\special{fp}%
%
\special{pn 13}%
\special{pa 4800 1200}%
\special{pa 1400 1400}%
\special{fp}%
%
\special{pn 13}%
\special{pa 3000 2000}%
\special{pa 4200 2600}%
\special{fp}%
%
\special{pn 13}%
\special{pa 2600 600}%
\special{pa 1400 1400}%
\special{dt 0.045}%
%
\special{pn 13}%
\special{pa 4800 1200}%
\special{pa 3000 2000}%
\special{dt 0.045}%
%
\special{pn 13}%
\special{pa 4200 2600}%
\special{pa 5600 1600}%
\special{pa 5600 800}%
\special{pa 3800 600}%
\special{pa 3800 600}%
\special{dt 0.045}%
\put(38.0000,-5.0000){\makebox(0,0)[lb]{$p_0$}}%
\put(26.0000,-5.0000){\makebox(0,0)[lb]{$p_1$}}%
\put(13.5000,-13.0000){\makebox(0,0)[lb]{$p_2$}}%
\put(42.0000,-13.8000){\makebox(0,0)[lb]{$p_3$}}%
\put(28.0000,-21.0000){\makebox(0,0)[lb]{$p_4$}}%
\put(43.0000,-27.0000){\makebox(0,0)[lb]{$p_5$}}%
%
\special{pn 4}%
\special{pa 4050 630}%
\special{pa 3400 1280}%
\special{fp}%
\special{pa 3990 630}%
\special{pa 3340 1280}%
\special{fp}%
\special{pa 3940 620}%
\special{pa 3280 1280}%
\special{fp}%
\special{pa 3880 620}%
\special{pa 3210 1290}%
\special{fp}%
\special{pa 3800 640}%
\special{pa 3150 1290}%
\special{fp}%
\special{pa 3760 620}%
\special{pa 3080 1300}%
\special{fp}%
\special{pa 3720 600}%
\special{pa 3020 1300}%
\special{fp}%
\special{pa 3660 600}%
\special{pa 2960 1300}%
\special{fp}%
\special{pa 3600 600}%
\special{pa 2890 1310}%
\special{fp}%
\special{pa 3540 600}%
\special{pa 2830 1310}%
\special{fp}%
\special{pa 3480 600}%
\special{pa 2770 1310}%
\special{fp}%
\special{pa 3420 600}%
\special{pa 2700 1320}%
\special{fp}%
\special{pa 3360 600}%
\special{pa 2640 1320}%
\special{fp}%
\special{pa 3300 600}%
\special{pa 2570 1330}%
\special{fp}%
\special{pa 3240 600}%
\special{pa 2510 1330}%
\special{fp}%
\special{pa 3180 600}%
\special{pa 2450 1330}%
\special{fp}%
\special{pa 3120 600}%
\special{pa 2380 1340}%
\special{fp}%
\special{pa 3060 600}%
\special{pa 2320 1340}%
\special{fp}%
\special{pa 3000 600}%
\special{pa 2260 1340}%
\special{fp}%
\special{pa 2940 600}%
\special{pa 2190 1350}%
\special{fp}%
\special{pa 2880 600}%
\special{pa 2130 1350}%
\special{fp}%
\special{pa 2820 600}%
\special{pa 2060 1360}%
\special{fp}%
\special{pa 2760 600}%
\special{pa 2000 1360}%
\special{fp}%
\special{pa 2700 600}%
\special{pa 1940 1360}%
\special{fp}%
\special{pa 2600 640}%
\special{pa 1870 1370}%
\special{fp}%
\special{pa 2510 670}%
\special{pa 1810 1370}%
\special{fp}%
\special{pa 2330 790}%
\special{pa 1750 1370}%
\special{fp}%
\special{pa 2150 910}%
\special{pa 1680 1380}%
\special{fp}%
\special{pa 1970 1030}%
\special{pa 1620 1380}%
\special{fp}%
\special{pa 1790 1150}%
\special{pa 1550 1390}%
\special{fp}%
%
\special{pn 4}%
\special{pa 1610 1270}%
\special{pa 1490 1390}%
\special{fp}%
\special{pa 4100 640}%
\special{pa 3470 1270}%
\special{fp}%
\special{pa 4150 650}%
\special{pa 3530 1270}%
\special{fp}%
\special{pa 4210 650}%
\special{pa 3590 1270}%
\special{fp}%
\special{pa 4260 660}%
\special{pa 3660 1260}%
\special{fp}%
\special{pa 4320 660}%
\special{pa 3720 1260}%
\special{fp}%
\special{pa 4370 670}%
\special{pa 3790 1250}%
\special{fp}%
\special{pa 4420 680}%
\special{pa 3850 1250}%
\special{fp}%
\special{pa 4480 680}%
\special{pa 3910 1250}%
\special{fp}%
\special{pa 4530 690}%
\special{pa 3980 1240}%
\special{fp}%
\special{pa 4590 690}%
\special{pa 4040 1240}%
\special{fp}%
\special{pa 4640 700}%
\special{pa 4100 1240}%
\special{fp}%
\special{pa 4690 710}%
\special{pa 4170 1230}%
\special{fp}%
\special{pa 4750 710}%
\special{pa 4230 1230}%
\special{fp}%
\special{pa 4800 720}%
\special{pa 4300 1220}%
\special{fp}%
\special{pa 4860 720}%
\special{pa 4360 1220}%
\special{fp}%
\special{pa 4910 730}%
\special{pa 4420 1220}%
\special{fp}%
\special{pa 4960 740}%
\special{pa 4490 1210}%
\special{fp}%
\special{pa 5020 740}%
\special{pa 4550 1210}%
\special{fp}%
\special{pa 5070 750}%
\special{pa 4610 1210}%
\special{fp}%
\special{pa 5130 750}%
\special{pa 4680 1200}%
\special{fp}%
\special{pa 5180 760}%
\special{pa 4770 1170}%
\special{fp}%
\special{pa 5230 770}%
\special{pa 4830 1170}%
\special{fp}%
\special{pa 5290 770}%
\special{pa 4830 1230}%
\special{fp}%
\special{pa 5340 780}%
\special{pa 3750 2370}%
\special{fp}%
\special{pa 5400 780}%
\special{pa 3790 2390}%
\special{fp}%
\special{pa 5450 790}%
\special{pa 3830 2410}%
\special{fp}%
\special{pa 5500 800}%
\special{pa 3870 2430}%
\special{fp}%
\special{pa 5560 800}%
\special{pa 3910 2450}%
\special{fp}%
\special{pa 5590 830}%
\special{pa 3950 2470}%
\special{fp}%
%
\special{pn 4}%
\special{pa 5590 890}%
\special{pa 3990 2490}%
\special{fp}%
\special{pa 5590 950}%
\special{pa 4030 2510}%
\special{fp}%
\special{pa 5590 1010}%
\special{pa 4070 2530}%
\special{fp}%
\special{pa 5590 1070}%
\special{pa 4110 2550}%
\special{fp}%
\special{pa 5590 1130}%
\special{pa 4150 2570}%
\special{fp}%
\special{pa 5590 1190}%
\special{pa 4220 2560}%
\special{fp}%
\special{pa 5590 1250}%
\special{pa 4380 2460}%
\special{fp}%
\special{pa 5590 1310}%
\special{pa 4590 2310}%
\special{fp}%
\special{pa 5590 1370}%
\special{pa 4800 2160}%
\special{fp}%
\special{pa 5590 1430}%
\special{pa 5010 2010}%
\special{fp}%
\special{pa 5590 1490}%
\special{pa 5220 1860}%
\special{fp}%
\special{pa 5590 1550}%
\special{pa 5430 1710}%
\special{fp}%
\special{pa 4830 1230}%
\special{pa 3710 2350}%
\special{fp}%
\special{pa 4770 1230}%
\special{pa 3670 2330}%
\special{fp}%
\special{pa 4680 1260}%
\special{pa 3630 2310}%
\special{fp}%
\special{pa 4570 1310}%
\special{pa 3590 2290}%
\special{fp}%
\special{pa 4460 1360}%
\special{pa 3550 2270}%
\special{fp}%
\special{pa 4350 1410}%
\special{pa 3510 2250}%
\special{fp}%
\special{pa 4240 1460}%
\special{pa 3470 2230}%
\special{fp}%
\special{pa 4140 1500}%
\special{pa 3430 2210}%
\special{fp}%
\special{pa 4030 1550}%
\special{pa 3390 2190}%
\special{fp}%
\special{pa 3920 1600}%
\special{pa 3350 2170}%
\special{fp}%
\special{pa 3810 1650}%
\special{pa 3310 2150}%
\special{fp}%
\special{pa 3700 1700}%
\special{pa 3270 2130}%
\special{fp}%
\special{pa 3600 1740}%
\special{pa 3230 2110}%
\special{fp}%
\special{pa 3490 1790}%
\special{pa 3190 2090}%
\special{fp}%
\special{pa 3380 1840}%
\special{pa 3150 2070}%
\special{fp}%
\special{pa 3270 1890}%
\special{pa 3110 2050}%
\special{fp}%
\special{pa 3160 1940}%
\special{pa 3070 2030}%
\special{fp}%
\end{picture}%
}
\caption{$D\subset \C P^1$}
\label{disk}
\end{minipage}
\hspace{20pt}
\begin{minipage}{.46\linewidth}
\scalebox{0.65}[0.65]{
\unitlength 0.1in
\begin{picture}( 42.5000, 23.1200)( 13.5000,-26.4200)
%
\special{pn 8}%
\special{sh 0.600}%
\special{ar 3800 600 42 42  0.0000000 6.2831853}%
%
\special{pn 8}%
\special{sh 0.600}%
\special{ar 2600 600 42 42  0.0000000 6.2831853}%
%
\special{pn 8}%
\special{sh 0.600}%
\special{ar 1400 1400 42 42  0.0000000 6.2831853}%
%
\special{pn 8}%
\special{sh 0.600}%
\special{ar 3000 2000 42 42  0.0000000 6.2831853}%
%
\special{pn 8}%
\special{sh 0.600}%
\special{ar 4800 1200 42 42  0.0000000 6.2831853}%
%
\special{pn 8}%
\special{sh 0.600}%
\special{ar 4200 2600 42 42  0.0000000 6.2831853}%
%
\special{pn 13}%
\special{pa 2600 600}%
\special{pa 3800 600}%
\special{fp}%
%
\special{pn 13}%
\special{pa 4800 1200}%
\special{pa 1400 1400}%
\special{fp}%
%
\special{pn 13}%
\special{pa 3000 2000}%
\special{pa 4200 2600}%
\special{fp}%
%
\special{pn 13}%
\special{pa 2600 600}%
\special{pa 1400 1400}%
\special{dt 0.045}%
%
\special{pn 13}%
\special{pa 4800 1200}%
\special{pa 3000 2000}%
\special{dt 0.045}%
%
\special{pn 13}%
\special{pa 4200 2600}%
\special{pa 5600 1600}%
\special{pa 5600 800}%
\special{pa 3800 600}%
\special{pa 3800 600}%
\special{dt 0.045}%
\put(38.0000,-5.0000){\makebox(0,0)[lb]{$P_0$}}%
\put(26.0000,-5.0000){\makebox(0,0)[lb]{$P_1$}}%
\put(13.5000,-13.0000){\makebox(0,0)[lb]{$P_2$}}%
\put(49.0000,-13.0000){\makebox(0,0)[lb]{$P_3$}}%
\put(27.7000,-21.0000){\makebox(0,0)[lb]{$P_4$}}%
\put(43.0000,-27.0000){\makebox(0,0)[lb]{$P_5$}}%
%
\special{pn 8}%
\special{sh 0.600}%
\special{ar 3400 1000 42 42  0.0000000 6.2831853}%
%
\special{pn 8}%
\special{pa 3400 1000}%
\special{pa 3800 600}%
\special{fp}%
%
\special{pn 8}%
\special{pa 3400 1000}%
\special{pa 2600 600}%
\special{fp}%
%
\special{pn 8}%
\special{pa 3400 1000}%
\special{pa 1400 1400}%
\special{fp}%
%
\special{pn 8}%
\special{pa 3400 1000}%
\special{pa 4800 1200}%
\special{pa 4800 1200}%
\special{fp}%
%
\special{pn 8}%
\special{pa 3400 1000}%
\special{pa 5200 1100}%
\special{pa 5000 1500}%
\special{pa 3000 2000}%
\special{pa 3000 2000}%
\special{pa 3000 2000}%
\special{fp}%
%
\special{pn 8}%
\special{pa 3400 1000}%
\special{pa 5400 1000}%
\special{pa 5400 1600}%
\special{pa 4200 2600}%
\special{pa 4200 2600}%
\special{pa 4200 2600}%
\special{fp}%
\put(33.5000,-12.5000){\makebox(0,0)[lb]{$Q_0$}}%
\put(37.0000,-8.7000){\makebox(0,0)[lb]{$\widetilde{\gamma}_0$}}%
\put(28.0000,-9.2000){\makebox(0,0)[lb]{$\widetilde{\gamma}_1$}}%
\put(45.0000,-22.0000){\makebox(0,0)[lb]{$\widetilde{\gamma}_5$}}%
\put(20.0000,-12.0000){\makebox(0,0)[lb]{$\widetilde{\gamma}_2$}}%
\put(36.0000,-20.3000){\makebox(0,0)[lb]{$\widetilde{\gamma}_4$}}%
\put(37.0000,-12.4000){\makebox(0,0)[lb]{$\widetilde{\gamma}_3$}}%
\end{picture}%
}
\caption{$\Omega_0\subset C$}
\label{lifting disk}
\end{minipage}
\end{center}
\end{figure}

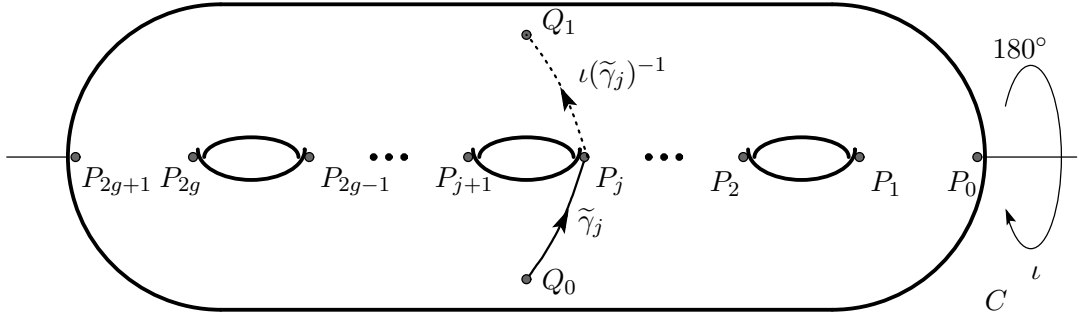
\begin{figure}[htbp]
\begin{center}
\unitlength 0.1in
\begin{picture}( 57.3000, 16.1000)(  5.5000,-26.0000)
%
\special{pn 13}%
\special{pa 3400 2440}%
\special{pa 3420 2414}%
\special{pa 3438 2388}%
\special{pa 3456 2362}%
\special{pa 3476 2336}%
\special{pa 3494 2310}%
\special{pa 3510 2282}%
\special{pa 3528 2256}%
\special{pa 3544 2228}%
\special{pa 3560 2200}%
\special{pa 3576 2172}%
\special{pa 3590 2144}%
\special{pa 3604 2116}%
\special{pa 3616 2086}%
\special{pa 3628 2058}%
\special{pa 3640 2028}%
\special{pa 3652 1998}%
\special{pa 3664 1968}%
\special{pa 3674 1938}%
\special{pa 3684 1906}%
\special{pa 3696 1876}%
\special{pa 3706 1846}%
\special{pa 3716 1816}%
\special{pa 3720 1800}%
\special{sp}%
%
\special{pn 20}%
\special{pa 1800 1000}%
\special{pa 5000 1000}%
\special{fp}%
\special{pa 5000 2600}%
\special{pa 1800 2600}%
\special{fp}%
%
\special{pn 20}%
\special{ar 1800 1800 800 800  1.5707963 4.7123890}%
%
\special{pn 20}%
\special{ar 5000 1800 800 800  4.7123890 6.2831853}%
\special{ar 5000 1800 800 800  0.0000000 1.5707963}%
%
\special{pn 20}%
\special{ar 1960 1760 280 160  6.2831853 6.2831853}%
\special{ar 1960 1760 280 160  0.0000000 3.1415927}%
%
\special{pn 20}%
\special{ar 1960 1800 248 104  3.1415927 6.2831853}%
%
\special{pn 20}%
\special{ar 3400 1760 280 160  6.2831853 6.2831853}%
\special{ar 3400 1760 280 160  0.0000000 3.1415927}%
%
\special{pn 20}%
\special{ar 4840 1760 280 160  6.2831853 6.2831853}%
\special{ar 4840 1760 280 160  0.0000000 3.1415927}%
%
\special{pn 20}%
\special{ar 3400 1800 248 104  3.1415927 6.2831853}%
%
\special{pn 20}%
\special{ar 4840 1800 248 104  3.1415927 6.2831853}%
%
\special{pn 8}%
\special{sh 0.600}%
\special{ar 3400 2440 24 24  0.0000000 6.2831853}%
%
\special{pn 8}%
\special{sh 0.600}%
\special{ar 1040 1800 24 24  0.0000000 6.2831853}%
%
\special{pn 8}%
\special{sh 0.600}%
\special{ar 5760 1800 24 24  0.0000000 6.2831853}%
%
\special{pn 8}%
\special{sh 0.600}%
\special{ar 5144 1800 24 24  0.0000000 6.2831853}%
%
\special{pn 8}%
\special{sh 0.600}%
\special{ar 3704 1800 24 24  0.0000000 6.2831853}%
%
\special{pn 8}%
\special{sh 0.600}%
\special{ar 2264 1800 24 24  0.0000000 6.2831853}%
%
\special{pn 8}%
\special{sh 0.600}%
\special{ar 1656 1800 24 24  0.0000000 6.2831853}%
%
\special{pn 8}%
\special{sh 0.600}%
\special{ar 3096 1800 24 24  0.0000000 6.2831853}%
%
\special{pn 8}%
\special{sh 0.600}%
\special{ar 4536 1800 24 24  0.0000000 6.2831853}%
%
\special{pn 8}%
\special{sh 0.600}%
\special{ar 3400 1160 24 24  0.0000000 6.2831853}%
%
\special{pn 20}%
\special{sh 1}%
\special{ar 2680 1800 10 10 0  6.28318530717959E+0000}%
\special{sh 1}%
\special{ar 2600 1800 10 10 0  6.28318530717959E+0000}%
\special{sh 1}%
\special{ar 2760 1800 10 10 0  6.28318530717959E+0000}%
%
\special{pn 20}%
\special{sh 1}%
\special{ar 4120 1800 10 10 0  6.28318530717959E+0000}%
\special{sh 1}%
\special{ar 4040 1800 10 10 0  6.28318530717959E+0000}%
\special{sh 1}%
\special{ar 4200 1800 10 10 0  6.28318530717959E+0000}%
%
\special{pn 8}%
\special{pa 1000 1800}%
\special{pa 680 1800}%
\special{fp}%
\special{pa 5800 1800}%
\special{pa 6280 1800}%
\special{fp}%
%
\special{pn 13}%
\special{pa 3720 1800}%
\special{pa 3710 1770}%
\special{pa 3700 1740}%
\special{pa 3690 1708}%
\special{pa 3680 1678}%
\special{pa 3668 1648}%
\special{pa 3658 1618}%
\special{pa 3646 1588}%
\special{pa 3634 1558}%
\special{pa 3622 1528}%
\special{pa 3610 1500}%
\special{pa 3596 1470}%
\special{pa 3582 1442}%
\special{pa 3568 1414}%
\special{pa 3552 1386}%
\special{pa 3536 1358}%
\special{pa 3520 1332}%
\special{pa 3502 1304}%
\special{pa 3484 1278}%
\special{pa 3466 1252}%
\special{pa 3448 1226}%
\special{pa 3428 1198}%
\special{pa 3410 1172}%
\special{pa 3400 1160}%
\special{sp -0.045}%
%
\special{pn 20}%
\special{pa 3600 2120}%
\special{pa 3608 2104}%
\special{fp}%
\special{sh 1}%
\special{pa 3608 2104}%
\special{pa 3560 2156}%
\special{pa 3584 2152}%
\special{pa 3596 2174}%
\special{pa 3608 2104}%
\special{fp}%
%
\special{pn 20}%
\special{pa 3600 1480}%
\special{pa 3592 1464}%
\special{fp}%
\special{sh 1}%
\special{pa 3592 1464}%
\special{pa 3604 1534}%
\special{pa 3616 1512}%
\special{pa 3640 1516}%
\special{pa 3592 1464}%
\special{fp}%
%
\special{pn 8}%
\special{ar 6040 1800 160 480  3.7310362 6.2831853}%
\special{ar 6040 1800 160 480  0.0000000 2.3561945}%
%
\special{pn 8}%
\special{pa 5928 2120}%
\special{pa 5912 2080}%
\special{fp}%
\special{sh 1}%
\special{pa 5912 2080}%
\special{pa 5918 2150}%
\special{pa 5932 2130}%
\special{pa 5956 2134}%
\special{pa 5912 2080}%
\special{fp}%
\put(34.8000,-11.6000){\makebox(0,0)[lb]{$Q_1$}}%
\put(34.8000,-25.2000){\makebox(0,0)[lb]{$Q_0$}}%
\put(56.0000,-20.0000){\makebox(0,0)[lb]{$P_0$}}%
\put(10.7000,-20.0000){\makebox(0,0)[lb]{$P_{2g+1}$}}%
\put(14.8000,-20.0000){\makebox(0,0)[lb]{$P_{2g}$}}%
\put(29.2000,-20.0000){\makebox(0,0)[lb]{$P_{j+1}$}}%
\put(43.6000,-20.0000){\makebox(0,0)[lb]{$P_2$}}%
\put(52.0000,-20.0000){\makebox(0,0)[lb]{$P_1$}}%
\put(37.6000,-20.0000){\makebox(0,0)[lb]{$P_j$}}%
\put(23.2000,-20.0000){\makebox(0,0)[lb]{$P_{2g-1}$}}%
\put(58.4000,-12.8000){\makebox(0,0)[lb]{$180^\circ$}}%
\put(60.4000,-24.4000){\makebox(0,0)[lb]{$\iota$}}%
\put(58.0000,-26.0000){\makebox(0,0)[lb]{$C$}}%
\put(36.7000,-22.2000){\makebox(0,0)[lb]{$\widetilde{\gamma}_j$}}%
\put(36.7000,-14.5000){\makebox(0,0)[lb]{$\iota(\widetilde{\gamma}_j)^{-1}$}}%
%
\special{pn 8}%
\special{pa 600 1800}%
\special{pa 550 1800}%
\special{ip}%
\end{picture}%
\caption{$e_j=\widetilde{\gamma}_j\cdot\iota(\widetilde{\gamma}_j)^{-1}$}
\label{hyperelliptic-curve}
\end{center}
\end{figure}

The hyperelliptic involution $\iota$
interchanges a point on $\Omega_0$
and the corresponding one on $\Omega_1$,
and fixes $P_i$, $i=0,1,\ldots,2g+1$.
We choose a base point $Q_0\in \Omega_0$
and denote $Q_1=\iota(Q_0)\in \Omega_1$.
Let $\gamma_j, j=0,1, \ldots, 2g+1$,
be a simple path in $D$ joining
$\pi(Q_0)$ and $p_j$.
We denote by $\widetilde{\gamma}_j$
the lift of $\gamma_j$ in $\Omega_0$ from $Q_0$ to $P_j$
(Figure \ref{lifting disk}, $g=2$).
Set $e_j=\widetilde{\gamma}_j\cdot\iota(\widetilde{\gamma}_j)^{-1}$,
where the product $\widetilde{\gamma}_j\cdot\iota(\widetilde{\gamma}_j)^{-1}$
indicates that we traverse $\widetilde{\gamma}_j$
first, then $\iota(\widetilde{\gamma}_j)^{-1}$.
It is a path in $C$ which is to be followed from 
$Q_0$ to $P_j$ and go to $Q_1$
in Figure \ref{hyperelliptic-curve}.
It is clear that
$e_{j_1}\cdot\iota(e_{j_2})$ is a loop
in $C$ at the base point $Q_0$.
Moreover we have the homotopy equivalences
relative to the base point $Q_0$
\begin{center}
$e_{j}\cdot\iota(e_{j})\sim 1$
\ and
\ $e_0 \cdot \iota(e_1)\cdot \cdots \cdot e_{2g} \cdot \iota(e_{2g+1}) \sim 1$.
\end{center}
We set
$a_i=e_{2i-1}\cdot\iota(e_{2i})$ and
$b_i=e_{2i-1}\cdot\iota(e_{2i-2})\cdot
\cdots\cdot e_1\cdot\iota(e_0)$,
and denote by $x_i$ and $y_i$
the homology classes of $a_i$ and $b_i$
respectively.
Then we have $\{x_i,y_i\}_{i=1,2,\ldots,g}$
is a symplectic basis of $H_1(C;\Z)$
in Figure 1 in \cite{T}.

\begin{figure}[htbp]
\begin{center}
\unitlength 0.1in
\begin{picture}( 52.3000,  5.2000)(  3.7000,-10.3000)
%
\special{pn 8}%
\special{sh 0.600}%
\special{ar 4000 800 46 46  0.0000000 6.2831853}%
%
\special{pn 8}%
\special{sh 0.600}%
\special{ar 5400 800 46 46  0.0000000 6.2831853}%
\put(35.7000,-10.4000){\makebox(0,0)[lb]{$\pi(Q_0)=\pi(Q_1)$}}%
\put(53.3000,-9.8000){\makebox(0,0)[lb]{$p_i$}}%
\put(44.8000,-6.8000){\makebox(0,0)[lb]{$\pi(e^\prime_i)$}}%
\put(56.0000,-12.0000){\makebox(0,0)[lb]{$\C P^1$}}%
%
\special{pn 13}%
\special{ar 5400 800 200 200  3.2904826 6.2831853}%
\special{ar 5400 800 200 200  0.0000000 2.9927027}%
%
\special{pn 13}%
\special{pa 4000 820}%
\special{pa 5200 820}%
\special{fp}%
\special{pa 4000 780}%
\special{pa 5200 780}%
\special{fp}%
%
\special{pn 20}%
\special{pa 4600 820}%
\special{pa 4610 820}%
\special{fp}%
\special{sh 1}%
\special{pa 4610 820}%
\special{pa 4544 800}%
\special{pa 4558 820}%
\special{pa 4544 840}%
\special{pa 4610 820}%
\special{fp}%
\special{pa 4760 780}%
\special{pa 4750 780}%
\special{fp}%
\special{sh 1}%
\special{pa 4750 780}%
\special{pa 4818 800}%
\special{pa 4804 780}%
\special{pa 4818 760}%
\special{pa 4750 780}%
\special{fp}%
%
\special{pn 20}%
\special{pa 5380 600}%
\special{pa 5370 600}%
\special{fp}%
\special{sh 1}%
\special{pa 5370 600}%
\special{pa 5438 620}%
\special{pa 5424 600}%
\special{pa 5438 580}%
\special{pa 5370 600}%
\special{fp}%
%
\special{pn 8}%
\special{sh 0.600}%
\special{ar 800 800 46 46  0.0000000 6.2831853}%
%
\special{pn 8}%
\special{sh 0.600}%
\special{ar 2200 800 46 46  0.0000000 6.2831853}%
\put(3.7000,-10.4000){\makebox(0,0)[lb]{$\pi(Q_0)=\pi(Q_1)$}}%
\put(21.3000,-9.8000){\makebox(0,0)[lb]{$p_i$}}%
\put(12.8000,-6.8000){\makebox(0,0)[lb]{$\pi(e_i)$}}%
\put(24.0000,-12.0000){\makebox(0,0)[lb]{$\C P^1$}}%
%
\special{pn 13}%
\special{pa 800 820}%
\special{pa 2200 820}%
\special{fp}%
\special{pa 800 780}%
\special{pa 2200 780}%
\special{fp}%
%
\special{pn 20}%
\special{pa 1400 820}%
\special{pa 1410 820}%
\special{fp}%
\special{sh 1}%
\special{pa 1410 820}%
\special{pa 1344 800}%
\special{pa 1358 820}%
\special{pa 1344 840}%
\special{pa 1410 820}%
\special{fp}%
\special{pa 1560 780}%
\special{pa 1550 780}%
\special{fp}%
\special{sh 1}%
\special{pa 1550 780}%
\special{pa 1618 800}%
\special{pa 1604 780}%
\special{pa 1618 760}%
\special{pa 1550 780}%
\special{fp}%
%
\special{pn 8}%
\special{pa 2800 800}%
\special{pa 3400 800}%
\special{fp}%
\special{pa 3300 700}%
\special{pa 3400 800}%
\special{fp}%
\special{pa 3400 800}%
\special{pa 3300 900}%
\special{fp}%
\end{picture}%
\caption{A deformation of $e_i$}
\label{deformation}
\end{center}
\end{figure}
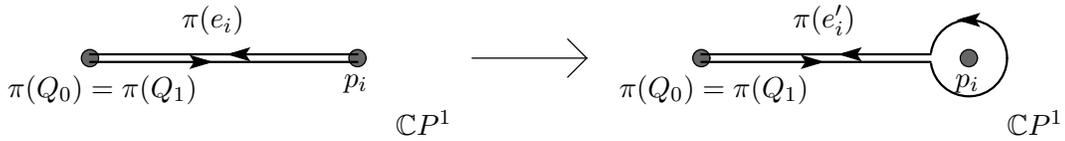

Let $H_{\Z_2}$ denote $H_1(C; \Z_2)$ and
$B$ branch locus $\{p_i\}_{i=0,1,\ldots,2g+1}$.
We deform the path $e_i$ in $C$ and denote it by
$e_i^{\prime}$ in $C\setminus \pi^{-1}(B)$ as follows.
The path $e_i^{\prime}$ avoids $P_i$
in a sufficiently small neighborhood
at $P_i$ so that 
$\pi(e_i^\prime)$ goes around $p_i$ 
which does not any $p_k$ with $k\neq i$
(Figure \ref{deformation}) and
the set of homology classes
$\{\pi(e_i^{\prime})\}_{i=0,1,\ldots, 2g}$
is a basis of $H_1(\C P^1\setminus B;\Z_2)$.
Moreover we have 
$\pi(e_0^{\prime})+\pi(e_1^{\prime})+\cdots+
\pi(e_{2g+1}^{\prime})=0\in H_1(\C P^1\setminus B;\Z_2)$.
Since the coefficients are in $\Z_2$, the homology class of
$e_i^\prime$ is independent of the choice of it.
From the $2$-sheeted covering $\pi$,
we have the well-defined homomorphism
$v_0\colon H\to H_1(\C P^1\setminus B;\Z_2)$ which factors through
$H_1(C\setminus \pi^{-1}(B);\Z)$ (Arnol'd \cite{A}).
We obtain the linear map
$v\colon H_{\Z_2}\to H_1(\C P^1\setminus B;\Z_2)$
induced naturally by $v_0$.
It immediately follows that
$v(x_i\ \mathrm{mod}\ 2)
=\pi(e_{2i-1}^{\prime})+\pi(e_{2i}^{\prime}),\
v(y_i\ \mathrm{mod}\ 2)
=\pi(e_0^{\prime})+\pi(e_1^{\prime})+\cdots+\pi(e_{2i-1}^{\prime})$,
and $v$ is injective.
The map $v$ gives the short exact sequence
$$
\xymatrix{
0 \ar[r] & H_{\Z_2} \ar[r]^{\hspace{-35pt}v} &
H_1(\C P^1\setminus B;\Z_2)
\ar[r] & \Z_2 \ar[r] & 0.
}
$$
Here the map $H_1(\C P^1\setminus B;\Z_2)\to \Z_2$ is the
augmentation map $\pi(e_i^\prime)\mapsto 1$.
Fix a Weierstrass point $P_{\nu}$.
Let $f_i$ denote $\pi(e_{\nu}^{\prime})+\pi(e_{i}^{\prime})$ for 
$i=0,1,\ldots, 2g+1$.
We remark that 
$f_i$ may be considered as an element of $H_{\Z_2}$
and $f_{\nu}=0$.
For $i=1,2,\ldots, g$, we have the identification
\begin{equation}\label{identification}
\left\{\begin{array}{l}
\displaystyle x_i=f_{2i-1}+ f_{2i} ,\\
\displaystyle y_i=f_0+f_1+\cdots +f_{2i-1},
\end{array}\right.
\end{equation}
in $H_{\Z_2}$.
It is clear that $f_0+f_1+\cdots +f_{2g+1}=0$.

For any hyperelliptic curve $C$
and Weierstrass point $P_{j}\in C$,
the hyperelliptic involution $\iota$
fixes $I_{P_{j}}$
and acts on $H_{\Z_2}$ as $(-1)$-times.
Then we have the value of $I_{P_{j}}$
is $0$ or $1/{2}\ \textnormal{mod}\ \Z$
from the equation $I_{P_j}\equiv (-1)^3I_{P_j}\ \mathrm{mod}\ \Z$.
We may consider $I_{P_{j}}\in \Hom((K\otimes H)_{\Z_2},\Z_2)$,
where $(K\otimes H)_{\Z_2}$ denotes $(K\otimes H)\otimes \Z_2$.

\section{Pointed harmonic volumes of hyperelliptic curves
and the moduli space of compact Riemann surfaces}
\label{Pointed harmonic volumes of hyperelliptic curves
and the moduli space of compact Riemann surfaces}
We recall some results about the moduli space
of compact Riemann surfaces.
Let $\Sigma_g$ be a closed oriented surface of genus $g$.
Its mapping class group, denoted here by $\Gamma^s_g$,
is the group of isotopy classes of orientation
preserving diffeomorphisms of $\Sigma_g$
which fix $s$ points on $\Sigma_g$ for $s=0,1$.
We denote $\Gamma_g=\Gamma^0_g$.
The group $\Gamma^1_g$ acts on the Teichm\"uller
space $\mathcal{T}^1_g$ of $\Sigma_g$ with a marked point
and the quotient space $\mathcal{M}^1_g$
is the moduli space of Riemann surfaces of genus $g$
with a marked point.
The group $\Gamma^1_g$ acts naturally on the first homology
group $H_1(\Sigma_g; \Z)$ of $\Sigma_g$.

Let $\mathcal{H}^1_g\subset \mathcal{M}^1_g$ be the moduli space of
hyperelliptic curves of genus $g$
with a marked Weierstrass point $P_{\nu}$.
For the rest of this paper, we suppose
that a marked point is a Weierstrass point.
The hyperelliptic mapping class group $\Delta^1_g$
is the subgroup of $\Gamma_g$ defined by
$$\{\varphi\in \Gamma_g; \varphi\iota=\iota\varphi,
\ \varphi(P_{\nu})=P_{\nu}\},$$
where $\iota$ is the hyperelliptic involution of $\Sigma_g$.
We have $\Delta^1_g\subset \Gamma^1_g$.
The moduli space $\mathcal{H}^1_g$
is known to be connected and
has a natural structure of a quasi-projective orbifold.
The group
$\Delta^1_g$
can be considered as its orbifold fundamental group.
For any $\Z\Delta^1_g$-module $M$,
we may consider the dual $M^{\ast}=\Hom(M,\Z_2)$
as a $\Z_2\Delta^1_g$-module in a natural way.
We denote $I_{\nu}=I_{P_{\nu}}$.
\begin{prop}\label{key fact}
We have
$$I_{\nu}\in H^0(\Delta^1_g; (K\otimes H)^{\ast}),$$
i.e. $I_{\nu}$ is a $\Delta^1_g$-invariant in the
dual $(K\otimes H)^{\ast}$.
\end{prop}
\begin{proof}
Let $\mathcal{L}$ be a locally constant sheaf with
a stalk $\Hom_{\Z}(K\otimes H, \Z_2)$.
In a similar way to Harris' method \cite{H-1},
$I_{\nu}$ varies in $\mathcal{H}^1_g$ continuously.
For any hyperelliptic curves,
$I_{\nu} \equiv 0$ or $1/{2}$ modulo $\Z$.
We remark that the pointed harmonic volume
is uniquely determined for any point on $\mathcal{H}^1_g$.
The locally constant sheaf
$\mathcal{L}$
has a global section
$\widetilde{I}_{\nu}$ associated to $I_{\nu}$.
Moreover $\mathcal{H}^1_g$ is arcwise connected.
Therefore $\widetilde{I}_{\nu}$
is a constant section of $\mathcal{L}$
and
$I_{\nu}$ is invariant under the action of the orbifold
fundamental group $\Delta_g^1$ of $\mathcal{H}^1_g$.
\end{proof}

\section{Pointed harmonic volumes of a hyperelliptic curve $C_0$}
\label{a hyperelliptic curve}
We compute the pointed harmonic volume
of a pointed hyperelliptic curve $(C_0,P_{\nu})$.
See \S 3 and 4 in \cite{T} for details.
We define the hyperelliptic curve $C_0$
by the equation $w^2=z^{2g+2}-1$.
We take $Q_i=(0, (-1)^{i}\sqrt{-1}), i=0, 1$,
and $P_j=(\zeta^j, 0), j=0,1, \ldots, 2g+1$,
where $\zeta=\exp\big(2\pi\sqrt{-1}/(2g+2)\big)$.
We define a path $e_j \colon[0,1] \to C_0, j=0, 1,\ldots ,2g+1$, by 
$$
\left\{
  \begin{array}{lll}
  \big( 2t\zeta^j,\sqrt{-1}\sqrt{1-(2t)^{2g+2}} \big)& \mathrm{for} &0\leq t\leq 1/{2}, \\
  \big( (2-2t)\zeta^j,-\sqrt{-1}\sqrt{1-(2-2t)^{2g+2}} \big)& \mathrm{for} & 1/{2}\leq t\leq 1.
  \end{array}
  \right.
$$

For $i=1,2, \ldots, g$,
we denote by $\omega_i$
a holomorphic $1$-form $z^{i-1}dz/{w}$ on $C_0$. 
It is known that $\{\omega_i\}_{i=1,2,\ldots ,g}$ is a basis
of the space of holomorphic $1$-forms on $C_0$.
Let $B(u,v)$ denote the beta function 
$\displaystyle \int_0^1x^{u-1}(1-x)^{v-1}dx$ for $u,v>0$.
For the normalization, we set 
$\omega^\prime_i=\displaystyle 
{(2g+2)\sqrt{-1}\over{2B(i/(2g+2),1/2)}}\omega_i$.
Then we have
$$
\int_{a_j}\omega^\prime_i=\zeta^{i(2j-1)}(1-\zeta^i)
\quad \textnormal{and}\ 
\int_{b_j}\omega^\prime_i={{\zeta^{2ij}-1}\over{\zeta^i}+1},
$$
where $i, j \in \{1,2, \ldots, g\}$.
The integral $\displaystyle\int_{\gamma} \omega^\prime_i$
depends only on the homology class of $\gamma$,
since $\omega^\prime_i$ is a closed $1$-form.

We compute
the iterated integrals of real harmonic
$1$-forms of $C_0$ with integral periods.
 Let $\Omega_a$ and $\Omega_b$ be the non-singular matrices
whose $(i,j)$-entries are
$$
\int_{a_j}\omega^\prime_i\quad
\mathrm{and}\ 
\int_{b_j}\omega^\prime_i
$$
respectively. 
We define real harmonic
$1$-forms $\alpha_i$ and $\beta_i, i=1,2, \ldots, g$, by
$$
\left(\begin{array}{c}
\alpha_1\\
 \vdots\\
\alpha_g
\end{array}
\right)=
\Re
\left(\begin{array}{c}
(\Omega_b)^{-1}
\left(\begin{array}{c}
\omega_1^\prime\\
 \vdots\\
 \omega_g^\prime
\end{array}
\right)
\end{array}
\right)
\hspace{5pt} \mathrm{and}\hspace{5pt}
\left(\begin{array}{c}
\beta_1\\
 \vdots\\
\beta_g 
\end{array}
\right)=-\Re
\left(\begin{array}{c}
(\Omega_a)^{-1}
\left(\begin{array}{c}
\omega_1^\prime\\
 \vdots\\
 \omega_g^\prime
\end{array}
\right)\end{array}
\right)
$$
respectively.
It is clear that
$\displaystyle\int_{a_j}\alpha_i=\int_{b_j}\beta_i=0$
and
$\displaystyle\int_{b_j}\alpha_i=\delta_{ij}=-\int_{a_j}\beta_i$.
Let $\Theta\colon H_1(C_0; \Z)\to H^1(C_0; \Z)$ denote
 the Poincar\'e dual.
We have
$\Theta(x_i)=\alpha_i$ and $\Theta(y_i)=\beta_i$ for
$i=1,2, \ldots, g$.
Hence, $\{\alpha_i, \beta_i\}_{i=1,2, \ldots, g}$
is a symplectic basis of $H^1(C_0; \Z)$.

Let $t_u$ be a complex number
$\displaystyle \sum_{p=1}^g\zeta^{up}$ for any integer $u$.
It is obvious that
$$
t_u=
\left\{
  \begin{array}{lll}
  g & \textnormal{for} &u\in (2g+2)\Z,\\
  -1 & \textnormal{for} &u\in 2\Z\setminus (2g+2)\Z,\\
  \displaystyle{{1+\zeta^u}\over{1-\zeta^u}} &\textnormal{for} & u\in 2\Z+1.
  \end{array}
  \right.
$$
Furthermore, $t_u$ is pure imaginary and
$t_{-u}=-t_u$ when $u$ is odd.
In addition to the formulas (1),(2),(3) and (4)
of Lemma 3.8 in \cite{T},
it is to show
\begin{lem}
\label{a hyperelliptic iterated integral}
On the curve $C_0$, we have
\begin{align*}
(5)\displaystyle \int_{a_k}\alpha_i\beta_j&
=\ \displaystyle\frac{-1}{2(g+1)^2}t_{2k-2i}(t_{2k-2j}-t_{2k}), \\
(6)\displaystyle \int_{b_k}\alpha_i\beta_j&
=\ \displaystyle\frac{-1}{2(g+1)^2}\sum_{u=1}^{k}\biggl\{
(t_{2u-2i-2}-t_{2u-2i})\sum_{v=1}^{j}t_{2v+2u-2j-2}\biggr\}.
\end{align*}
Here $i, j, k \in \{1,2, \ldots, g\}$.
\end{lem}
\begin{proof}
We compute the case (5) in the following way.
Let $A_{j,m}$ and $B_{i,l}$ be
$(j,m)$ and $(i,l)$-entries of
$(\Omega_a)^{-1}$ and $(\Omega_b)^{-1}$
respectively.
\begin{align*}
&\int_{a_k}\alpha_i\beta_j=
\int_{a_k}-\Re\biggl
(\sum_{l=1}^{g}B_{i,l}\omega^{\prime}_l\biggr)
\Re\biggl(\sum_{m=1}^{g}A_{j,m}\omega^{\prime}_m\biggr)\\
=&
-\frac{1}{4}
\int_{a_k}
\sum_{l,m=1}^{g}
\biggl(
B_{i,l}A_{j,m}\omega^{\prime}_l\omega^{\prime}_m
+B_{i,l}\overline{A}_{j,m}\omega^{\prime}_l\overline{\omega}^{\prime}_m
+\overline{B}_{i,l}A_{j,m}\overline{\omega}^{\prime}_l\omega^{\prime}_m
+\overline{B}_{i,l}\overline{A}_{j,m}
\overline{\omega}^{\prime}_l\overline{\omega}^{\prime}_m
\biggr)\\
=&
-\frac{1}{2}
\Re\Biggl\{
\sum_{l,m=1}^{g}
\biggl(
B_{i,l}A_{j,m}\int_{a_k}\omega^{\prime}_l\omega^{\prime}_m
+B_{i,l}\overline{A}_{j,m}
\int_{a_k}\omega^{\prime}_l\overline{\omega}^{\prime}_m
\biggr)
\Biggr\}.
\end{align*}
We use Lemma 3.5 in \cite{T}
and calculate
\begin{align*}
&(g+1)^2\sum_{l,m=1}^{g}
B_{i,l}A_{j,m}\int_{a_k}\omega^{\prime}_l\omega^{\prime}_m\\
=&
\sum_{l,m=1}^{g}
\zeta^{-2il}(1+\zeta^l)
\frac{\zeta^m(-1+\zeta^{-2jm})}{1-\zeta^m}
{1\over2}\zeta^{(l+m)(2k-1)}(1-2\zeta^m+\zeta^{l+m})\\
=&\frac{1}{2}
\sum_{m=1}^{g}
\frac{1-\zeta^{2jm}}{1-\zeta^m}
\zeta^{m(2k-2j)}
\sum_{l=1}^{g}\zeta^{l(2k-2i-1)}
(1+\zeta^l)(1-\zeta^m-\zeta^m(1-\zeta^l))\\
=&\frac{1}{2}
\sum_{m=1}^{g}
\sum_{v=2k-2j}^{2k-1}\zeta^{mv}
\Big\{(1-\zeta^m)(t_{2k-2i-1}+t_{2k-2i})
-\zeta^m(t_{2k-2i-1}-t_{2k-2i+1})\Big\}\\
=&\frac{1}{2}
\sum_{m=1}^{g}\biggl\{
\frac{1-\zeta^{2jm}}{1-\zeta^m}
\zeta^{m(2k-2j)}(1-\zeta^m)(t_{2k-2i-1}+t_{2k-2i})
-\sum_{v=2k-2j}^{2k-1}\zeta^{m(v+1)}
(t_{2k-2i-1}-t_{2k-2i+1})
\biggr\}\\
=&\frac{1}{2}
\biggl\{(t_{2k-2i-1}+t_{2k-2i})(t_{2k-2j}-t_{2k})
-(t_{2k-2i-1}-t_{2k-2i+1})\sum_{v=2k-2j}^{2k-1}
t_{v+1}
\biggr\}\\
=&\frac{1}{2}
\biggl\{(t_{2k-2i-1}+t_{2k-2i})(t_{2k-2j}-t_{2k})
-(t_{2k-2i-1}-t_{2k-2i+1})\sum_{v=2k-2j+1}^{2k}
t_{v}
\biggr\}.
\end{align*}
Similarly, we have
\begin{align*}
&(g+1)^2\sum_{l,m=1}^{g}
B_{i,l}\overline{A}_{j,m}
\int_{a_k}\omega^{\prime}_l\overline{\omega}^{\prime}_m\\
=&\frac{1}{2}
\biggl\{(t_{2k-2i-1}+t_{2k-2i})(t_{2k-2j}-t_{2k})
-(t_{2k-2i-1}-t_{2k-2i+1})\sum_{v=2k-2j+1}^{2k}
t_{-v}
\biggr\}.
\end{align*}
Therefore, we obtain the result
\begin{align*}
&\int_{a_k}\alpha_i\beta_j\\
=&\frac{-1}{2(g+1)^2}
\frac{1}{2}\Re\biggl\{
2(t_{2k-2i-1}+t_{2k-2i})(t_{2k-2j}-t_{2k})
-2(t_{2k-2i-1}-t_{2k-2i+1})
\sum_{\begin{subarray}{c}v=2k-2j+1\\ 
\mathrm{even}\end{subarray}}^{2k}t_{v}
\biggr\}\\
=&\frac{-1}{2(g+1)^2}
t_{2k-2i}(t_{2k-2j}-t_{2k}).
\end{align*}
Similarly we compute the case (6).
\end{proof}

Using the symplectic basis
$\{x_i,y_i\}_{i=1,2,\cdots,g}\subset H_1(C;\Z)$
stated in \S \ref{Hyperelliptic curves},
we choose a basis of $K$ as follows:
$$\left\{
\begin{array}{ccc}
(1) & z_i\otimes z^{\prime}_j & (i\neq j)\\
(2) & x_i\otimes y_i-x_1\otimes y_1 & (i\neq 1)\\
(3) & x_i\otimes y_i+y_i\otimes x_i & (i=1,2,\ldots, g)\\
(4) & z_i\otimes z_i & (i=1,2,\ldots, g)
\end{array}
\right\},
$$
where $z_i$ denotes $x_i$ or $y_i$ and so on.
By the definition of the pointed harmonic volume
$I_{\nu}$,
we obtain
$$I_{\nu}((x_i\otimes y_i+y_i\otimes x_i)\otimes z^{\prime\prime}_k)
\equiv 0\ \mathrm{mod}\ \Z\quad \textnormal{for any } i,k,$$
and
$$
I_{\nu}(z_i\otimes z_i\otimes z^{\prime\prime}_k)\equiv 
\left\{
\begin{array}{cl}
1/{2}\ \mathrm{mod}\ \Z & \textnormal{if}\ z_i\otimes z_i\otimes 
z^{\prime\prime}_k=
x_i\otimes x_i\otimes y_i\ \textnormal{or}\ y_i\otimes y_i\otimes x_i,\\
0\ \mathrm{mod}\ \Z & \textnormal{otherwise},
\end{array}
\right.
$$
for any hyperelliptic curve $C$.
It is enough to consider the case (1) and (2).
For the rest of this paper,
we omit $\mathrm{mod}\ \Z$,
unless otherwise stated.

We compute the pointed harmonic volume
of $(C_0,Q_0)$.
From Lemma \ref{a hyperelliptic iterated integral},
Lemma 3.8 in \cite{T}
and the equation $\displaystyle\int_{e_j}\eta=0$
(Lemma 4.2 in \cite{T}), it is to show
\begin{prop}\label{p.h.v. with a pointed h.c.}
Case (1).
If $i\neq k$ and $j\neq k$, then we have
$$
I_{Q_0}(z_i\otimes z^{\prime}_j\otimes z^{\prime\prime}_k)
\equiv 0.
$$
If $i= k$ or $j= k$, then we have
$$
\begin{array}{lcll}
I_{Q_0}(x_i\otimes x_j\otimes y_i)&
\equiv &\mu,&\\
I_{Q_0}(x_i\otimes y_j\otimes y_i)&
\equiv &\left\{
\begin{array}{l}
(g-j+1)\mu \\
(2g-j+2)\mu
\end{array}
\right.
&
\begin{array}{l}
\textnormal{if }i<j,\\
\textnormal{if }i>j,
\end{array}
\\
 I_{Q_0}(y_i\otimes x_j\otimes x_i)&
\equiv &(2g+1)\mu,&\\
 I_{Q_0}(y_i\otimes y_j\otimes x_i)&
\equiv &\left\{
\begin{array}{l}
(g+j+1)\mu \\
j\mu 
\end{array}
\right.
&
\begin{array}{l}
\textnormal{if }i<j,\\
\textnormal{if }i>j,
\end{array}
\end{array}
$$

Case (2).
If $i\neq k$ and $k\neq 1$, then we have
$$
I_{Q_0}((x_i\otimes y_i-x_1\otimes y_1)\otimes z^{\prime\prime}_k)
\equiv 0.
$$
If $i=k$ or $k=1$, then we have
$$
\begin{array}{lcl}
I_{Q_0}((x_i\otimes y_i-x_1\otimes y_1)\otimes x_i)&
\equiv &(g+2)\mu,\\
I_{Q_0}((x_i\otimes y_i-x_1\otimes y_1)\otimes y_i)&
\equiv &(2g-i+2)\mu,\\
I_{Q_0}((x_i\otimes y_i-x_1\otimes y_1)\otimes x_1)&
\equiv &g\mu,\\
I_{Q_0}((x_i\otimes y_i-x_1\otimes y_1)\otimes y_i)&
\equiv &(g+2)\mu.\\
\end{array}
$$
Here we denote $\mu=1/{(2g+2)}$.
\end{prop}
\begin{rem}
From Remark \ref{invarinant},
we do not need to compute
$I_{Q_0}(x_j\otimes x_i\otimes y_i)$,
$I_{Q_0}((y_i\otimes x_i-y_1\otimes x_1)\otimes x_i)$
and so on.
\end{rem}

We calculate the difference between
$I_{\nu}$ and $I_{Q_0}$.
For $h_1\otimes h_2\otimes h_3\in K\otimes H$,
we set $\Lambda_{\nu}(h_1\otimes h_2\otimes h_3)
:= I_{\nu}(h_1\otimes h_2\otimes h_3)-
I_{Q_0}(h_1\otimes h_2\otimes h_3)\ \textnormal{mod}\ \Z$.
Let $\ell_{\nu}\colon [0,1]\to C_0$
be a path
$t\mapsto (t\zeta^{\nu}, \sqrt{-1}\sqrt{1-t^{2g+2}})\in C_0$.
It is clear that
$\ell_{\nu}^{-1}\cdot e_j\cdot \ell_{\nu}$'s
are loops in $C_0$ at the base point $P_{\nu}$.
From the equation (2.2) in \cite{H-1},
we have
\begin{lem}\label{difference}
$$\Lambda_{\nu}(h_1\otimes h_2\otimes h_3)
\equiv (h_1,h_3)\int_{\ell_{\nu}}h_2-(h_2,h_3)\int_{\ell_{\nu}}h_1
\quad \mathrm{mod}\ \Z.$$
\end{lem}
It is clear that
$$\int_{\ell_{\nu}}\alpha_i=\frac{1}{2(g+1)}
\Re(t_{\nu-2i}+t_{\nu-2i+1})\quad \textnormal{and}\ 
\int_{\ell_{\nu}}\beta_i=\frac{-1}{2(g+1)}
\Re\Biggl(\sum_{u=\nu-2i+1}^{\nu}t_u\Biggr).$$
These equations and Lemma \ref{difference} give the following Lemma.
\begin{lem}\label{difference2}
Case (1).
If $i\neq k$ and $j\neq k$, then we have
$$
\Lambda_{\nu}(z_i\otimes z^{\prime}_j\otimes z^{\prime\prime}_k)
\equiv 0.
$$
If $i= k$ or $j=k$, then we have
$$
\begin{array}{lcll}
\Lambda_{\nu}(x_i\otimes x_j\otimes y_i)&
\equiv &\left\{
\begin{array}{l}
g\mu\\
(2g+1)\mu
\end{array}
\right.
&
\begin{array}{l}
\textnormal{if }\nu=1\ \textnormal{or}\ 2,\\
\textnormal{if }\nu\neq 1\ \textnormal{and}\ 2,
\end{array}
\\
\Lambda_{\nu}(x_i\otimes y_j\otimes y_i)&
\equiv &\left\{
\begin{array}{l}
j\mu\\
(g+j+1)\mu
\end{array}
\right.
&
\begin{array}{l}
\textnormal{if }\nu>2i-1,\\
\textnormal{if }\nu\leq 2i-1,
\end{array}
\\
\Lambda_{\nu}(y_i\otimes x_j\otimes x_i)&
\equiv &\left\{
\begin{array}{l}
(g+2)\mu\\
\mu 
\end{array}
\right.
&
\begin{array}{l}
\textnormal{if }\nu=1\ \textnormal{or}\ 2,\\
\textnormal{if }\nu\neq 1\ \textnormal{and}\ 2,
\end{array}
\\
\Lambda_{\nu}(y_i\otimes y_j\otimes x_i)&
\equiv &\left\{
\begin{array}{l}
(2g-j+2)\mu\\
(g-j+1)\mu
\end{array}
\right.
&
\begin{array}{l}
\textnormal{if }\nu>2i-1,\\
\textnormal{if }\nu\leq 2i-1.
\end{array}
\end{array}
$$
Case (2).
If $i\neq k$ and $k\neq 1$, then we have
$$
\Lambda_{\nu}((x_i\otimes y_i-x_1\otimes y_1)\otimes z^{\prime\prime}_k)
\equiv 0.
$$
If $i= k$ or $k=1$, then we have
$$
\begin{array}{lcll}
\Lambda_{\nu}((x_i\otimes y_i-x_1\otimes y_1)\otimes x_i)&
\equiv &\left\{
\begin{array}{l}
g\mu\\
(2g+1)\mu
\end{array}
\right.
&
\begin{array}{l}
\textnormal{if }\nu=2j-1\ \textnormal{or}\ 2j,\\
\textnormal{if }\nu\neq 2j-1\ \textnormal{and}\ 2j,
\end{array}
\\
\Lambda_{\nu}((x_i\otimes y_i-x_1\otimes y_1)\otimes y_i)&
\equiv &\left\{
\begin{array}{l}
i\mu\\
(g+i+1)\mu
\end{array}
\right.
&
\begin{array}{l}
\textnormal{if }\nu>2i-1,\\
\textnormal{if }\nu\leq 2i-1,
\end{array}
\\
\Lambda_{\nu}((x_i\otimes y_i-x_1\otimes y_1)\otimes x_1)&
\equiv &\left\{
\begin{array}{l}
(g+2)\mu\\
\mu
\end{array}
\right.
&
\begin{array}{l}
\textnormal{if }\nu=1\ \textnormal{or}\ 2,\\
\textnormal{if }\nu\neq 1\ \textnormal{and}\ 2,
\end{array}
\\
\Lambda_{\nu}((x_i\otimes y_i-x_1\otimes y_1)\otimes y_1)&
\equiv &\left\{
\begin{array}{l}
(2g+1)\mu\\
g\mu
\end{array}
\right.
&
\begin{array}{l}
\textnormal{if }\nu>1,\\
\textnormal{if }\nu\leq 1.
\end{array}
\end{array}
$$
\end{lem}
By combining Proposition \ref{p.h.v. with a pointed h.c.}
and Lemma \ref{difference2},
we have the pointed harmonic volume $I_{\nu}$
of $(C_0,P_{\nu})$.
\begin{thm}\label{a pointed harmonic volume}
Case (1).
Elements of $K\otimes H$ at which the value of
the pointed harmonic volumes $I_{\nu}$
are $1/{2}\ \textnormal{mod}\ \Z$ are given by
$$
\begin{array}{lll}
x_i\otimes x_j\otimes y_i, & x_j\otimes x_i\otimes y_i
& \textnormal{if }\nu=2j-1\ \textnormal{or}\ 2j,\\
x_i\otimes y_j\otimes y_i, & y_j\otimes x_i\otimes y_i
& \textnormal{if }(i<j,\nu>2j-1)\ \textnormal{or}\ (i>j,\nu \leq 2j-1),\\
y_i\otimes x_j\otimes x_i, & x_j\otimes y_i\otimes x_i
& \textnormal{if }\nu=2j-1\ \textnormal{or}\ 2j,\\
y_i\otimes y_j\otimes x_i, & y_j\otimes y_i\otimes x_i
& \textnormal{if }(i<j,\nu>2j-1)\ \textnormal{or}\ (i>j,\nu \leq 2j-1).
\end{array}
$$
The values at the other elements are $0\ \textnormal{mod}\ \Z$.

Case (2).
Elements of $K\otimes H$ at which the value of
the pointed harmonic volumes $I_{\nu}$
are $1/{2}\ \textnormal{mod}\ \Z$ are given by
$$
\begin{array}{lll}
(x_i\otimes y_i-x_1\otimes y_1)\otimes x_i,
 & (y_i\otimes x_i-y_1\otimes x_1)\otimes x_i
& \textnormal{if }\nu\neq 2i-1\ \textnormal{and}\ 2i,\\
(x_i\otimes y_i-x_1\otimes y_1)\otimes y_i,
 & (y_i\otimes x_i-y_1\otimes x_1)\otimes y_i
& \textnormal{if }\nu\leq 2i-1,\\
(x_i\otimes y_i-x_1\otimes y_1)\otimes x_1,
 & (y_i\otimes x_i-y_1\otimes x_1)\otimes x_1
& \textnormal{if }\nu\neq 1\ \textnormal{and}\ 2,\\
(x_i\otimes y_i-x_1\otimes y_1)\otimes y_1,
& (y_i\otimes x_i-y_1\otimes x_1)\otimes y_1
& \textnormal{if }\nu>1.
\end{array}
$$
The values at the other elements are $0\ \textnormal{mod}\ \Z$.
\end{thm}

From Proposition \ref{key fact},
this theorem can be extended to any hyperelliptic curve $C$
with Weierstrass base points.
But this extension is complicated.
We reconsider Theorem \ref{a pointed harmonic volume}
from a combinatorial viewpoint.
We apply an element $A\in K\otimes H$
to the identification (\ref{identification}) in the group
$(K\otimes H)_{\Z_2}$.
Then we have
$(A\ \mathrm{mod}\ 2)=
\sum_{p,q,r\neq\nu}A_{p,q,r}f_p\otimes f_q\otimes f_r$,
where $A_{p,q,r}\in \Z_2=\{0,1\}$.
The notation $\sharp$ means the cardinality of a set.
A counting function
$\kappa_{\nu}\colon K\otimes H\to \dfrac{1}{2}\Z/{\Z}=\{0,1/{2}\}$
is well-defined by
$$\kappa_{\nu}(A):=\dfrac{1}{2}
\Big(\sharp\{(p,q,r); A_{p,q,r}=1,
\sharp\{p,q,r\}=2\}\Big)\ \mathrm{mod}\ \Z.$$
Here $\sharp\{p,q,r\}=2$
means
$p=q\neq r$ or $q=r\neq p$ or $r=p\neq q$.
By the long but easy computation,
we obtain the correspondence.
\begin{cor}\label{cor}
On the curve $C_0$, we have
$$
I_{\nu}(A)\equiv \kappa_{\nu}(A)\ \mathrm{mod}\ \Z.
$$
\end{cor}
\begin{ex}

\begin{enumerate}
\item If $A=x_i\otimes x_j\otimes y_i$ $(i<j$ and $\nu=2j-1)$,
we have
\begin{align*}
\kappa(A)
=&\ \kappa((f_{2i-1}+f_{2i})
\otimes f_{2j}\otimes 
(f_0+f_1+\cdots+f_{2i-1}))\\
\equiv&\ \kappa(\ft{2i-1}{2j}{2i-1})=1/{2}.
\end{align*}
\item If $A=x_i\otimes x_j\otimes y_i$ $(i>j$ and $2i<\nu)$,
we have
\begin{align*}
\kappa(A)
=&\ \kappa((f_{2i-1}+f_{2i})\otimes (f_{2j-1}+f_{2j})\otimes 
(f_0+f_1+\cdots+f_{2i-1}))\\
\equiv&\ \kappa(\ft{2i-1}{2j-1}{2j-1}+\ft{2i-1}{2j-1}{2i-1}
+\ft{2i-1}{2j}{2j}\\
&\qquad+\ft{2i-1}{2j}{2i-1}+\ft{2i}{2j-1}{2j-1}+\ft{2i}{2j}{2j})\\
=&\ 
1/{2}+1/{2}+1/{2}+1/{2}+1/{2}+1/{2}
\equiv 0.
\end{align*}
\end{enumerate}
\end{ex}

\section{A combinatorial formula of $I_{\nu}$}
\label{formula}
In this section, we compute the pointed harmonic volume
$I_{\nu}=I_{P_\nu}$ of $(C,P_{\nu})$ by another combinatorial way.
Let $S_{2g+1}$ be the $(2g+1)$-th symmetric group.
Using the natural projection
$\Delta^1_g\to S_{2g+1}$,
the group $H_{\Z_2}$ is naturally considered
as a $\Z_2S_{2g+1}$-module (Arnol'd, V. I. \cite{A}).
From the slight modification of
Lemma 5.5 and Proposition 5.7 in \cite{T},
we have
\begin{lem}\label{nontrivial element}
$$H^0(\Delta^1_g; (K\otimes H)^{\ast})=
H^0(S_{2g+1}; (H^{\otimes 3})^{\ast})=\Z_2.$$
Moreover the unique nontrivial element
$\psi_{\nu}\in H^0(S_{2g+1}; (H^{\otimes 3})^\ast)$
is an $S_{2g+1}$-homomorphism
$H^{\otimes 3}\to \Z_2$ defined by
$$\psi_{\nu}(f_i\otimes f_j\otimes f_k)=
\left\{
\begin{array}{ll}
1 &if\ \sharp\{i,j,k\}=2,\\
0 &otherwise,
\end{array}
\right.$$
for any $i,j,k$ except for $\nu$.
\end{lem}
From Lemma \ref{nontrivial element}, we have
\begin{thm}\label{a combinatorial formula}
For $A\in K\otimes H$, we have
$$
I_{\nu}(A)\equiv \kappa_{\nu}(A)\ \mathrm{mod}\ \Z.
$$
\end{thm}
Using the equation
$f_i=\pi(e^{\prime}_{\nu})+\pi(e^{\prime}_{i})$,
we obtain
$A=\sum_{p,q,r}A^{\prime}_{p,q,r}\pi(e^{\prime}_{p})
\otimes \pi(e^{\prime}_{q})\otimes \pi(e^{\prime}_{r})$.
Another counting function
$\kappa_{\nu}^{\prime}\colon K\otimes H\to
\dfrac{1}{2}\Z/{\Z}=\{0,1/{2}\}$
is defined by
$$\kappa_{\nu}^{\prime}(A):=\dfrac{1}{2}\Big(
\sharp\{(p,q,r); A^{\prime}_{p,q,r}=1,
\sharp\{p,q,r\}=2,\ p,q,r\neq\nu\}\Big)\ \mathrm{mod}\ \Z.$$
\begin{cor}
$$
I_{\nu}(A)\equiv\kappa_{\nu}^{\prime}(A)\ \mathrm{mod}\ \Z.
$$
\end{cor}
\begin{proof}
We use the notation
$e(p,q,r)=\pi(e^{\prime}_{p})
\otimes \pi(e^{\prime}_{q})\otimes \pi(e^{\prime}_{r})$
only here.
The equation
\begin{align*}
f_p\otimes f_q\otimes f_r=&
e(p,q,r)+e(p,q,\nu)+e(p,\nu,r)+e(p,\nu,\nu)\\
&+e(\nu,q,r)+e(p,q,\nu)+e(\nu,\nu,r)+e(\nu,\nu,\nu)
\end{align*}
gives $\kappa_{\nu}(A)\equiv\kappa_{\nu}^{\prime}(A)$.
\end{proof}

\end{document}